\documentclass[12pt,reqno]{amsart}

\usepackage{amsmath}
\usepackage{amsfonts}
\usepackage{amssymb}
\usepackage{mathrsfs}

\usepackage{a4}

\usepackage{xy}
\xyoption{all}

\usepackage{enumerate}

\usepackage{exscale,txfonts}

\renewcommand{\implies}{\Rightarrow}
\renewcommand{\iff}{\Leftrightarrow}
\newcommand{\too}{\longrightarrow}
\newcommand{\mapstoo}{\longmapsto}
\newcommand{\ox}{\otimes}
\renewcommand{\bar}{\overline}

\DeclareMathOperator{\Int}{Int}
\DeclareMathOperator{\Sym}{Sym} 
\DeclareMathOperator{\End}{End} 
\DeclareMathOperator{\Trd}{Trd}

\newcommand{\op}{{\mathrm{op}}}
\newcommand{\<}{\langle}
\renewcommand{\>}{\rangle}
\newcommand{\x}{\times}
\DeclareMathOperator{\Skew}{Skew}

\DeclareMathOperator{\Hom}{Hom}
\DeclareMathOperator{\Nil}{Nil}
\DeclareMathOperator{\dotcup}{\dot\cup}
\DeclareMathOperator{\Herm}{\mathscr{H}\!}
\DeclareMathOperator{\type}{\mathrm{type}}
\DeclareMathOperator*{\bigcupdot}{\dot\bigcup}

\newcommand{\wh}[1]{\ensuremath{\widehat{#1}}}
\newcommand{\wwh}{\raisebox{-.5ex}{\wh{\phantom{x}}\,}}

\newcommand{\R}{\mathbb{R}}

\newcommand{\N}{\mathbb{N}}
\newcommand{\Z}{\mathbb{Z}}

\newcommand{\M}{\mathscr{M}}
\newcommand{\CN}{\mathscr{N}}

\newcommand{\ve}{\varepsilon}
\newcommand{\vf}{\varphi}
\newcommand{\vt}{\vartheta}

\newcommand{\s}{\sigma}

\newcommand{\GL}{\mathrm{GL}}

\newcommand{\bbar}{-}

\newcommand{\ad}{\mathrm{ad}}
\newcommand{\id}{\mathrm{id}}

\newcommand{\sm}{\setminus}

\renewcommand{\epsilon}{\varepsilon}

\DeclareMathOperator{\sign}{sign}

\newcommand{\nc}{{\M}}

\newcommand{\simtoo}{\overset\sim\longrightarrow}

\newcommand{\trs}{\mathrm{Tr}^*}
\newcommand{\Tr}{\mathrm{Tr}}

\newtheorem{lemma}{Lemma}[section]
\newtheorem{theorem}[lemma]{Theorem}
\newtheorem{proposition}[lemma]{Proposition}
\newtheorem{corollary}[lemma]{Corollary}

\theoremstyle{definition}
\newtheorem{definition}[lemma]{Definition}

\newtheorem{remark}[lemma]{Remark}

\newtheorem{example}[lemma]{Example}

\title{Signatures of hermitian forms and the Knebusch Trace Formula}
\author{Vincent Astier}
\author{Thomas Unger}
\address{School of Mathematical Sciences\\ University College Dublin\\ Belfield\\
Dublin~4\\ Ireland} 
\email{vincent.astier@ucd.ie, thomas.unger@ucd.ie}


\begin{document}

\begin{abstract}  
Signatures of quadratic forms have been generalized to hermitian forms over algebras with involution. 
In the literature this is done via Morita theory, which causes sign ambiguities in certain cases.
In this paper, a hermitian version of the Knebusch Trace Formula is established and used as a main tool 
to resolve these ambiguities.
\end{abstract}

\maketitle

\section{Introduction}

In this paper we study signatures of hermitian forms over central simple algebras with involution
of any kind, defined over formally real fields. These 
generalize  the classical signatures of quadratic forms.

Following \cite{BP2} we do this via extension to real closures and Morita equivalence. 
This leads to the notion of $M$-signature of  hermitian forms in Section~\ref{sec:general}.
We study its properties, make a detailed analysis of the impact of choosing different real closures and different Morita equivalences  and show in particular that sign changes can occur. This motivates the search for a more intrinsic notion of signature, where such sign changes do not occur. 

In Section~\ref{sec:H} we define such a signature, the $H$-signature, which only depends on the choice of a tuple of hermitian forms,
mimicking   the fact that in quadratic form theory the form $\<1\>$ always has positive signature.
The $H$-signature  generalizes the definition of signature in \cite{BP2} and is in particular well-defined when  the involution becomes hyperbolic after scalar extension to a real closure of the base field, addressing an issue with the definition proposed in \cite{BP2}.

Our main tool is a generalization of the Knebusch trace formula to  $M$-signatures of hermitian forms,
which we establish in Section~\ref{ktf-M}. In Section~\ref{sec:cont} we show that the 
total $H$-signature of 
a hermitian form is a continuous map and in Section~\ref{sec:ktf-H} we prove the Knebusch trace formula for $H$-signatures.

\section{Preliminaries}

\subsection{Algebras with Involution} The general reference for this section is  \cite[Chapter~I]{BOI}. 

Let $F$ be a field of characteristic different from $2$.
An \emph{$F$-algebra with involution} is a pair $(A,\s)$ consisting of a finite-dimen\-sional   $F$-algebra $A$ with centre  $Z(A)$ and an
$F$-linear map $\s:A\to A$ such that $\s(xy)=\s(y)\s(x)$ for all $x,y\in A$ and $\s^2=\id_A$. The involution $\s$ is 
either of \emph{the first kind} or of   \emph{the second kind}. In the first case
$A$ is simple, $Z(A)=F$ and $\s|_F=\id_F$. In the second case there are two possibilities: either
$A$ is simple,  $Z(A)=K=F(\sqrt{d})$ for some $d\in F^\x$  and  $\s|_{K}$ is the nontrivial $F$-automorphism of $K$, or $(A,\s)\simeq (B\x B^\op, \wwh)$ with $B$ a simple $F$-algebra, $Z(A)\simeq F\x F$  a double-field and $\wwh$  the exchange involution, defined by 
$\wh{(x,y^\op)}=(y,x^\op)$ for all $x,y\in B$. We call $(A,\s)$ \emph{degenerate} if it is isomorphic to $(B\x B^\op, \wwh)$.

Consider the $F$-subspaces
$\Sym(A,\s)=\{a\in A\mid \s(a)=a\}$  and $\Skew(A,\s)=\{a\in A\mid \s(a)=-a\}$
of $A$. 
Then $A=\Sym(A,\s) \oplus \Skew(A,\s)$.
Assume that $\s$ is of the first kind. Then $\dim_F(A)=m^2$ for some positive integer $m$. Furthermore,
$\s$ is either \emph{orthogonal} (or, \emph{of type $+1$}) if $\dim_F\Sym(A,\s)=m(m+1)/2$, or  \emph{symplectic} (or, \emph{of  type $-1$}) if $\dim_F\Sym(A,\s)=m(m-1)/2$.  If $\s$ is of the second kind, then $\dim_F(A) = 2m^2$ for some positive integer $m$ and 
 $\dim_F\Sym(A,\s)=\dim_F\Skew(A,\s)=m^2$. Involutions of the second kind are also called \emph{unitary}.

Let $\s$ and $\tau$ be two involutions on $A$ that have the same restriction to $Z(A)$.
By the Skolem-Noether theorem they differ by an inner automorphism: 
$\tau = \Int(u)\circ \s$
for some $u\in A^\x$, uniquely determined up to a factor in $F^\x$,
 such that $\s(u)=u$ if $\s$ and $\tau$ are both orthogonal, both symplectic or both unitary and  $\s(u)=-u$ if one of $\s$, $\tau$ is orthogonal and the other symplectic. Here $\Int(u)(x):=uxu^{-1}$ for $x\in A$.

\subsection{$\ve$-Hermitian Spaces and Forms} 
The general references for this section are \cite[Chapter I]{Knus} and \cite[Chapter~7]{Sch}, both for rings with involution. Treatments of the central simple and division cases can also be found in \cite{GB} and \cite{Le}, respectively.

Let $(A,\s)$ be an $F$-algebra with involution. Let $\ve\in \{-1,1\}$. An \emph{$\ve$-hermitian space} over $(A,\s)$ is a pair $(M,h)$, where $M$ is a finitely generated right $A$-module (which is automatically projective since $A$ is semisimple) and $h:M\x M\too A$ is a sesquilinear form such that $h(y,x)=\ve \s(h(x,y))$ for all $x,y \in M$.
We call $(M,h)$ a \emph{hermitian space} when $\ve=1$ and a \emph{skew-hermitian space} when $\ve=-1$. If $(A,\s)$ is a field equipped with the identity map, we say \emph{\textup{(}skew-\textup{)} symmetric bilinear} space instead of (skew-) hermitian space.

Consider the left $A$-module $M^*=\Hom_A(M,A)$ as a right $A$-module via the involution $\s$. The form $h$ induces an $A$-linear map $h^*:M\to M^*, x\mapsto h(x,\cdot)$. We call $(M,h)$ \emph{nonsingular} if $h^*$ is an isomorphism. 
All spaces occurring in this paper are assumed to be nonsingular.  We often simply write $h$ instead of $(M,h)$ and speak of a \emph{form} instead of a space. 

By Wedderburn theory, $M$ decomposes into a direct sum of $k$ simple right $A$-modules, for some $k\in \N$ which we call the \emph{rank} of $h$.

If $A=D$ is a division algebra (so that $M\simeq D^n$ for some integer $n$)  such that $(D,\s,\ve)\not=(F,\id_F, -1)$, then $h$ can be diagonalized: there exist invertible elements $a_1,\ldots, a_n\in \Sym(D,\s)$ such that, after a change of basis,
$h(x,y)=\sum_{i=1}^n \s(x_i) a_i y_i$, for all $x=(x_1,\ldots,x_n), y=(y_1,\ldots, y_n) \in D^n$.
In this case we use the shorthand notation
$h=\<a_1,\ldots, a_n\>_\s$,
which resembles the usual notation for diagonal quadratic forms. 

If $A$ is not a division algebra  we can certainly consider diagonal hermitian forms 
$h=\<a_1,\ldots, a_n\>_\s$
defined on the  free $A$-module $A^n$, where  $n\in \N$  and  
$a_1,\ldots, a_n$ are invertible elements in $\Sym(A,\s)$. 
We call $n$ the \emph{dimension} of $h$.
Note that not all hermitian forms over $(A,\s)$  are diagonalizable. 

Witt cancellation and Witt decomposition hold for  $\ve$-hermitian spaces 
over 
$(A,\s)$.
Let $W_\ve(A,\s)$ denote the Witt group 
of Witt classes of $\ve$-hermitian spaces over $(A,\s)$.
When $\ve=1$ we drop the subscript and simply write  $W(A,\s)$. 
We denote the usual Witt ring of $F$ by $W(F)$. We find it convenient to
identify forms over $(A,\s)$ with their classes in $W_\ve(A,\s)$.

\begin{lemma} \label{lem:witt-isom}
Let $(A,\s)$ be an $F$-algebra with involution.
\begin{enumerate}[$(i)$]
\item If $\s$ is of the first kind and $\dim_F(A)$ is odd, then $W(A,\s)\simeq W(F)$ and $W_{-1}(A,\s)=0$. 
\item If $\s$ is of the first kind and  $\dim_F(A)$ is even, then $W_{-1}(A,\s)\simeq W(A,\tau)$ for some involution $\tau$ of opposite type to $\s$.
\item If $\s$ is of the second kind, then $W_{-1} (A,\s) \simeq W(A,\s)$.
\item If $(A,\s)$ is degenerate, then $W(A,\s)=0$.
\end{enumerate}
\end{lemma}

\begin{proof} $(i)$ It follows from the assumptions that $A\simeq M_n(F)$ and that $\s$ is necessarily orthogonal. By hermitian Morita theory
(see Section~\ref{sec:hmt}) we have
that $W(A,\s)\simeq W(F)$ and $W_{-1}(A,\s) \simeq W_{-1}(F,\id_F)$. Since skew-symmetric forms over $F$ are hyperbolic we have $W_{-1}(F,\id_F)=0$.

$(ii)$ By \cite[3.A]{BOI} there exists an involution $\tau$  of opposite type to $\s$ since  $A$ has involutions of both types under our assumptions. As observed earlier we then have $\tau\circ\s = \Int(u)$ for some $u\in A^\x$ with $\s(u)=-u$.  
Let $h$ be a skew-hermitian form over $(A,\s)$. Then $uh$ is a hermitian form over $(A,\tau)$. The one-to-one correspondence $h\mapsto uh$ respects isometries, orthogonal sums and hyperbolicity and so induces the indicated isomorphism.

$(iii)$ Let $u\in Z(A)$, $u\not= 0$ be such that $\s(u)=-u$. For example, let $u=\sqrt{d}$ if $Z(A)=F(\sqrt{d})$ and let $u=(-1,1)$ if $Z(A)\simeq F\x F$. 
The one-to-one correspondence $h\mapsto uh$ induces the indicated isomorphism.

$(iv)$ Assume that $(A,\s) \simeq (B\x B^\op, \wwh)$ and let $h:M\x M \to B\x B^\op$ be hermitian with respect to $\wwh$. Let $e_1=(1,0)$ and $e_2=(0,1)$. Then $M=M_1\oplus M_2$ with $M_i:=Me_i$ $(i=1,2)$ and  $h$ is hyperbolic since $h|_{M_1\x M_1}=0$ and 
$M_1^\perp = M_1$ (cf. \cite[Chapter~I, Corollary~3.7.3]{Knus}), which can be verified by direct computation.\qed
\end{proof}

In view of Lemma~\ref{lem:witt-isom}$(iv)$ degenerate $F$-algebras with involution are not 
interesting in the context of this paper. Therefore 
we redefine \emph{$F$-algebra with involution} to mean  non-degenerate $F$-algebra with involution. Observe that such an algebra with involution may become degenerate over a field extension of $F$.

\subsection{Adjoint Involutions}\label{adj.invo} The general reference for this section is  \cite[4.A]{BOI}.

Let $(A,\s)$ be an $F$-algebra with involution.
 Let $(M,h)$ be an $\ve$-hermitian space over $(A,\s)$. The algebra $\End_A(M)$ is again simple with centre $Z(A)$ since $M$ is finitely generated \cite[1.10]{BOI}.
The involution $\ad_h$ on $\End_A(M)$, defined by
$h(x, f(y))=h(\ad_h(f)(x),y)$, for all $x,y \in M$,  and all $f\in \End_A(M)$
is called the \emph{adjoint involution} of $h$. The involutions $\s$ and $\ad_h$ are of the same kind and $\s(\alpha)=\ad_h(\alpha)$ for all $\alpha\in Z(A)$. In case $\ad_h$ and $\s$ are of the first kind we also have
$\type(\ad_h)=\ve\, \type(\s)$.
Furthermore, every involution on $\End_A(M)$ is of the form $\ad_h$ for some $\ve$-hermitian form $h$ over $(A,\s)$ and the correspondence between $\ad_h$ and $h$ is unique up to a multiplicative factor  in $F^\x$ in the sense that $\ad_h=\ad_{\lambda h}$ for every $\lambda \in F^\x$.

Let $(A,\s)$ be an $F$-algebra with involution.
By a theorem of Wedderburn there exists a division algebra $D$ (unique up to isomorphism) with centre $Z(A)$ and a  finite-dimensional
right $D$-vector space $V$ such that $A \simeq \End_D(V)$.
Thus $A \simeq M_m(D)$   for some positive integer $m$. Furthermore there  exists an involution $\vt$ on $D$ of the same kind as $\s$ and an $\ve_0$-hermitian   form $\vf_0$ over $(D,\vt)$  with $\ve_0\in \{-1,1\}$
such that $(A,\s)$ and $(\End_D(V), \ad_{\vf_0})$ are isomorphic as algebras with involution. In matrix form
$\ad_{\vf_0}$ is described as follows:
$\ad_{\vf_0}(X)=\Phi_0\vt^t(X)\Phi_0^{-1}$, for all $X\in M_m(D)$,
where $\vt^t(X):=(\vt(x_{ij}))^t$ for $X=(x_{ij})$ and
$\Phi_0\in \GL_m(D)$ is the Gram matrix of $\vf_0$. Thus $\vt^t(\Phi_0)=\ve_0 \Phi_0$. 

\subsection{Hermitian Morita Theory}\label{sec:hmt}
We refer to \cite[\S1]{BP1}, \cite{FM},  \cite[Chapters 2--3]{GB}, \cite[Chapter~I, \S9]{Knus},  or \cite{Le1} for more details. 

Let $(M,h)$ be an $\ve$-hermitian space over $(A,\s)$.
One can show that the algebras with involution $(\End_A(M),\ad_h)$ and $(A, \s)$ are Morita equivalent:  for every $\mu\in \{-1,1\}$
there is an equivalence between the categories
$\Herm_\mu(\End_A(M),\ad_h)$ and $\Herm_{\ve \mu}(A,\s)$ of non-singular $\mu$-hermitian forms over $(\End_A(M),\ad_h)$   and non-singular $\ve \mu$-hermitian forms over $(A,\s)$, respectively (where the morphisms are given by isometry), cf. \cite[Chapter~I, Theorem~9.3.5]{Knus}. This equivalence  respects isometries, orthogonal sums and hyperbolic forms. It
 induces  a group isomorphism
$W_\mu (\End_A(M), \ad_h)  \simeq  W_{\ve\mu}(A,\s)$.
The
Morita equivalence and the isomorphism are 
not canonical.  

The algebras with involution $(A,\s)$ and $(D, \vt)$ are also Morita equivalent. 
An example of such a Morita equivalence is obtained by composing the following three non-canonical equivalences of categories, the last two of which we will call \emph{scaling} and \emph{collapsing}. For computational purposes  we describe them in matrix form. We follow the approach of \cite{LU2}:
\begin{small}
\begin{equation}\label{comp:morita}
\Herm_{\ve} (A,\s)\xrightarrow{\phantom{xxx}}
\Herm_{\ve} (M_m(D), \ad_{\vf_0}) \xrightarrow{\text{ scaling  }} 
\Herm_{\ve_0 \ve} (M_m(D), \vt^t) \xrightarrow{\text{ collapsing  }}
 \Herm_{\ve_0 \ve} (D, \vt).
\end{equation}
\end{small}

\textbf{Scaling:} Let $(M,h)$ be an $\ve$-hermitian space over $(M_m(D), \ad_{\vf_0})$. Scaling is given by 
$(M,h)\mapstoo (M, \Phi_0^{-1}h)$.
Note that $\Phi_0^{-1}$ is only determined up to a scalar factor in $F^\x$ since 
$\ad_{\vf_0}=\ad_{\lambda \vf_0}$ for any $\lambda\in F^\x$ and that replacing $\Phi_0$ by $\lambda \Phi_0$ results in a different Morita equivalence.
\medskip

\textbf{Collapsing:} Recall that $M_m(D) \simeq \End_D(D^m)$ and
that we always have $M  \simeq (D^m)^k \simeq M_{k,\, m}(D)$ for some integer $k$. Let $h: M\x M\too M_m(D)$ be an $\ve_0 \ve$-hermitian form with respect to $\vt^t$. Then
$h(x,y)=\vt^t(x) B y$, for all $x,y\in M_{k,\, m}(D)$,
where $B\in M_k(D)$ satisfies $\vt^t(B)=\ve_0\ve B$, so that $B$ determines an $\ve_0\ve$-hermitian form $b$ over $(D,\vt)$. Collapsing is then given by 
$(M,h)\mapstoo (D^k, b)$.

\section{Signatures of Hermitian Forms}

\subsection{Signatures of Forms: the Real Closed Case}\label{sec:special}
Let $R$ be a real closed field, $C=R(\sqrt{-1})$ (which is algebraically closed) and $H=(-1,-1)_R$ Hamilton's quaternion division algebra over $R$. We recall the definitions of signature for the various types of forms (all assumed to be nonsingular) over $R$, $C$, $R\x R$ and $H$.  We will use them
in the definition of the $M$-signature of a hermitian form over $(A,\s)$ in Section~\ref{sec:general}.

\begin{enumerate}[(a)]
\item Let $b$ be a symmetric bilinear (or quadratic) form over $R$. Then $b\simeq \< \alpha_1, \ldots, \alpha_n\>$ for some $n\in \N$ and $\alpha_i \in \{-1,1\}$. We let 
$\sign b:= \sum_{i=1}^n \alpha_i$. 
By Sylvester's Law of Inertia, $\sign b$ is well-defined. 

\item Let $b$ be a skew-symmetric form over $R$. Then $b$ is hyperbolic and we let
$\sign b:=0$.

\item Let $h$ be a hermitian form over $(C,-)$, where $\bar{\sqrt{-1}}=-\sqrt{-1}$. Then $h\simeq \< \alpha_1, \ldots, \alpha_n\>_{-}$ for some $n\in \N$ and $\alpha_i \in \{-1,1\}$. By a theorem of Jacobson \cite{J}, $h$ is up to isometry uniquely determined by the symmetric bilinear form $b_h:=2\x  \< \alpha_1, \ldots, \alpha_n\>$ defined over $R$. We let
$\sign h:=  \frac{1}{2} \sign b_h = \sign \< \alpha_1, \ldots, \alpha_n\>$.

\item Let $h$ be a hermitian form over $(R\x R, \wwh)$, where $\wh{(x,y)}=(y,x)$ is the exchange involution. Then $h$ is hyperbolic by Lemma~\ref{lem:witt-isom}$(iv)$ and we let
$\sign h:=0$.

\item Let $h$ be a hermitian form over $(H,-)$, where $-$ denotes quaternion conjugation. Then $h\simeq \< \alpha_1, \ldots, \alpha_n\>_{-}$ for some $n\in \N$ and $\alpha_i \in \{-1,1\}$. By a theorem of Jacobson \cite{J}, $h$ is up to isometry uniquely determined by the symmetric bilinear form 
$b_h:=4\x  \< \alpha_1, \ldots, \alpha_n\>$ defined over $R$. We let
$\sign h:=  \frac{1}{4} \sign b_h = \sign \< \alpha_1, \ldots, \alpha_n\>$.

\item Let $h$ be a skew-hermitian form over $(H,-)$, where $-$ denotes quaternion conjugation. Then $h$ is a torsion form \cite[Chapter~10, Theorem~3.7]{Sch} and we let
$\sign h:=0$.
\end{enumerate}
Note that the cases of skew-hermitian forms over $(C,-)$ and $(R\x R, \wwh)$ can be reduced to (c) and (d), respectively, by Lemma~\ref{lem:witt-isom}.

\begin{remark}\label{Z-iso} 
The  signature maps defined in (a), (c) and (e)  above give rise to unique group isomorphisms
$W(R)\simeq \Z,\ W(C,-) \simeq \Z$, and $W(H,-)\simeq \Z$
such that $\sign \<1\>=1$, $\sign \<1\>_{-}=1$ and $\sign \<1\>_{-}=1$, respectively. 
In addition, we have the group isomorphisms
$W_{-1}(R,\id_R)= W(R\x R, \wwh)=0$,  $W_{-1}(H,\bbar)\simeq \Z/2\Z$.
See also  \cite[Chapter~I, 10.5]{Knus}. 
\end{remark}

\subsection{The $M$-Signature of a Hermitian Form}\label{sec:general}
Our approach in this section is inspired by  \cite[\S3.3,\S3.4]{BP2}. 
We only consider hermitian forms over $(A,\s)$, cf. Lemma~\ref{lem:witt-isom}. 

Let $F$
be a formally real field and let $(A,\s)$ be an $F$-algebra with involution. 
Consider an ordering $P\in X_F$, the space of orderings of $F$.
By a \emph{real closure of $F$ at $P$} we mean a field embedding $\iota : F
\rightarrow K$, where $K$ is real closed, $\iota(P) \subseteq K^2$ and $K$ is
algebraic over $\iota(F)$.

Let
$h$ be a hermitian form over $(A,\s)$.   
Choose a real closure $\iota : F
\rightarrow F_P$ of $F$ at $P$, and use it to extend scalars from $F$ to $F_P$:
\[W(A,\s) \too W(A\ox_F F_P, \s\ox \id),\ h\mapstoo h\ox F_P:= (\id_A \ox \iota)^*(h)  ,\]
where the tensor product is along $\iota$.
The extended algebra with involution $(A\ox_F F_P, \s\ox \id)$ is  Morita equivalent to an algebra with involution $(D_P, \vt_P)$, chosen as follows:

\begin{enumerate}[$(i)$]
\item If $\s$ is of the first kind, $D_P$ is equal to one of $F_P$ or $H_P:=(-1,-1)_{F_P}$. Furthermore, we may choose $(D_P, \vt_P)= (F_P, \id_{F_P})$ in the first case and 
$(D_P,\vt_P)= (H_P,-) $ in the second case by Morita theory (scaling).

\item If $\s$ is of the second kind,  recall that $Z(A)=F(\sqrt{d})$. Now if $d<_P0$, then $D_P$ is equal to $F_P(\sqrt{-1})$, whereas  if $d>_P 0$, then $D_P$ is equal to $F_P\x F_P$ and $A\ox_F F_P$ is a direct product of two simple algebras. Furthermore, we may choose 
$(D_P, \vt_P) = (F_P(\sqrt{-1}),   -)$ in the first case and $(D_P, \vt_P) = 
(F_P\x F_P, \wwh)$ in the second case, again by Morita theory (scaling).
\end{enumerate}
Note that $\vt_P$ is of the same kind as $\s$ in each case.

The extended involution $\s\ox \id_{F_P}$ is adjoint to
an $\ve_P$-hermitian form $\vf_P$ over $(D_P, \vt_P)$  where $\ve_P=-1$ if one of $\s$ and $\vt_P$ is orthogonal and the other is symplectic, whereas $\ve_P=1$ if $\s$ and $\vt_P$ are of the same type, i.e. both orthogonal, symplectic or unitary. 

Now choose any Morita equivalence
\begin{equation}\label{eq:morita}
\M: \Herm(A\ox_F F_P, \s\ox \id) \too \Herm_{\ve_P}(D_P, \vt_P)
\end{equation}
with $(D_P,\vt_P) \in \{ (F_P, \id_{F_P}), (H_P,-), (F_P(\sqrt{-1}),   -), (F_P\x F_P, \wwh)\}$, which exists by the analysis above.
This Morita equivalence induces an isomorphism, which we again denote by $\M$, 
\begin{equation}\label{mor.iso}
\M: W(A\ox_F F_P, \s\ox \id) \simtoo W_{\ve_P}(D_P, \vt_P).
\end{equation}

\begin{definition} \label{def:sign} 
Let $P \in X_F$. Fix a real closure $\iota : F \rightarrow F_P$ of $F$ at $P$
and a Morita equivalence $\M$ as above. Define the \emph{$M$-signature of $h$
at $(\iota,\M)$}, denoted
$\sign_\iota^\nc h$, as follows:
\[\sign_\iota^\nc h:= \sign \M(h\ox F_P),\]
where $\sign \M(h\ox F_P)$ can be computed as shown in Section~\ref{sec:special}.
\end{definition}

This definition relies on two choices: firstly the choice of the real closure $\iota: F \rightarrow
F_P$ of $F$ at $P$ and secondly the choice of the Morita equivalence $\M$. 
Note that there is no canonical choice for $\M$. 
We now study the dependence of the $M$-signature on the choice of $\iota$ and
$\M$.

Let $\iota_1 : F \rightarrow L_1$ and $\iota_2 : F \rightarrow L_2$ be two real
closures of $F$ at $P$, and let $(D_1,\vt_1)$ and $\ve_1$ play the role of 
$(D_P,\vt_P)$ and $\ve_P$, respectively, obtained above when $\iota$ is replaced by $\iota_1$. Let 
$\M_1: \Herm(A \ox_F L_1, \s \ox \id) \too \Herm_{\ve_1}(D_1, \vt_1)$
be a fixed Morita equivalence. 
By the Artin-Schreier theorem \cite[Chapter~3,Theorem~2.1]{Sch} there is a unique isomorphism
$\rho : L_1 \rightarrow L_2$ such that $\rho \circ
\iota_1 = \iota_2$. It extends to an isomorphism $\id \ox \rho : (A \ox_F L_1, \s \ox \id)
\rightarrow (A \ox_F L_2, \s \ox \id)$. The isomorphism $\rho$ also extends canonically to
$D_1 \in\{L_1, (-1,-1)_{L_1}, L_1(\sqrt{-1}), L_1\x L_1 \}$. Consider the $L_2$-algebra with 
involution $(D_2,\vt_2):=(\rho(D_1), \rho\circ \vt_1\circ \rho^{-1})$. We define $\rho(\M_1)$ to be the Morita equivalence 
from $(A \ox_F L_2, \s \ox \id)$ to $(D_2, \vt_2)$,  described by the following  diagram: 
\begin{equation*}
\xymatrix{
\Herm(A\ox_F L_1, \s\ox\id) \ar[rr]^--{\M_1} \ar[d]^--{(\id\ox \rho)^*}& & \Herm_{\ve_1}(D_1,\vt_1)\ar[d]^--{\rho^*}\\
\Herm(A\ox_F L_2, \s\ox\id) \ar[rr]^--{\rho(\M_1)} &  &\Herm_{\ve_1}(D_2,\vt_2) \\
}
\end{equation*}
 
\begin{proposition}[Change of Real Closure]\label{prop:iota} 
With notation as above we have  for every $h\in W(A,\s)$,
\[\sign \M_1(h\ox L_1) =
\sign \rho(\M_1)(h\ox L_2),\]
in other words \[\sign_{\iota_1}^{\nc_1} h =
\sign_{\iota_2}^{\rho(\nc_1)} h.\]
\end{proposition}

\begin{proof}  The statement is trivially true when $\ve_1=-1$, by cases (b), (d) and (f) in Section~\ref{sec:special}, so we may assume that $\ve_1=1$.
Consider the diagram
\begin{equation*}
\begin{aligned}[c] 
\xymatrix{
                &W(A\ox_F L_1, \s\ox\id) \ar[rr]^--{\M_1}_--{\sim} \ar[dd]^--{(\id\ox \rho)^*}& & W(D_1,\vt_1)\ar[dd]^--{\rho^*}  \ar[rd]^--{\sign}& \\
W(A,\s)\ar[ur]^--{(\id\ox\iota_1)^*} \ar[dr]_--{(\id\ox\iota_2)^*}   &         &        &    & \Z\\
               &W(A\ox_F L_2, \s\ox\id) \ar[rr]^--{\rho(\M_1)}_--{\sim} &  &W(D_2,\vt_2) \ar[ur]_--{\sign} &  \\
}
\end{aligned}
\end{equation*}
The left triangle commutes by the definition of $\rho$. The square commutes by the definition of $\rho(\M_1)$. The right triangle commutes since $\rho^*(\<1\>_{\vt_1}) = \<1\>_{\vt_2}$ and by Remark~\ref{Z-iso}.
The statement follows.\qed
\end{proof}

\begin{proposition}[Change of Morita Equivalence]
\label{prop:morita}  
Let $\M_1$ and $\M_2$ be two different Morita equivalences as in \eqref{eq:morita}. Then there exists $\delta \in \{-1,1\}$ such that for every $h \in W(A,\s)$,
\[ \sign_\iota^{\M_1} h = \delta  \sign_\iota^{\M_2} h.\]
\end{proposition}
  
\begin{proof}  Note that $\ve_1=\ve_2$. The statement is trivially true when $\ve_1=-1$,  by cases (b), (d) and (f) in Section~\ref{sec:special}, so we may assume that $\ve_1=1$.
The two different Morita equivalences give rise to two different group isomorphisms
\[m_i: W(A\ox_F F_P, \s\ox \id)  \simtoo \Z\qquad (i=1,2),\]
by \eqref{mor.iso} and Remark~\ref{Z-iso}. The map $m_1\circ m_2^{-1}$ is an automorphism of $\Z$ and is therefore equal to $\id_\Z$ or $-\id_\Z$.\qed
\end{proof}

Propositions~\ref{prop:iota} and \ref{prop:morita} immediately imply

\begin{corollary}\label{iota:morita} 
Let $\iota_1 : F \rightarrow L_1$ and $\iota_2 : F \rightarrow L_2$ be two real
closures of $F$ at $P$ and let $\M_1$ and $\M_2$ be two different Morita equivalences as in \eqref{eq:morita}. Then there exists $\delta \in \{-1,1\}$ such that for every $h \in W(A,\s)$,
\[ \sign_{\iota_1}^{\M_1} h = \delta  \sign_{\iota_2}^{\M_2} h.\]
\end{corollary}

The following result easily follows    from the properties of Morita equivalence:

\begin{proposition} \label{rem:morsign} \mbox{}
\begin{enumerate}[$(i)$]
\item Let $h$ be a hyperbolic form over $(A,\s)$, then
$\sign^\nc_\iota h=0$.

\item Let $h_1$ and
$h_2$ be hermitian forms over $(A,\sigma)$, then
$\sign^\nc_\iota (h_1\perp h_2)=\sign^\nc_\iota h_1 + \sign^\nc_\iota h_2$.

\item The $M$-signature at $(\iota, \M)$, $\sign^\nc_\iota$, induces a homomorphism of additive groups 
$W(A,\s)\too \Z$.

\item Let $h$ be a hermitian form over $(A,\s)$ and $q$ a quadratic form over $F$, then
$\sign_\iota^\M (q\ox h) = \sign_P q \cdot \sign_\iota^\M h$,
where $\sign_P q$ denotes the usual signature of the quadratic form $q$ at $P$.
\end{enumerate}
\end{proposition}

\begin{definition}
Let $h$ be a hermitian form over $(A,\s)$. From Definition~\ref{def:sign} and Section~\ref{sec:special} it follows that $\sign^\nc_\iota h$ is automatically zero whenever $P$ belongs to the  following subset of $X_F$, which we call set of \emph{nil-orderings}:
\[
\Nil[A,\s]:=
\begin{cases}
\{ P\in X_F \mid D_P= H_P\} & \text{if $\s$ is orthogonal}\\
\{ P\in X_F \mid D_P= F_P\} & \text{if $\s$ is symplectic}\\
\{ P\in X_F \mid D_P= F_P\x F_P\} & \text{if $\s$ is unitary}\\
\end{cases},
\]
where the square brackets indicate that $\Nil[A,\s]$ depends only on the Brauer class of $A$ and the type of $\s$.
\end{definition}

\subsection{The $H$-Signature of a Hermitian Form}\label{sec:H}
It follows from Corollary~\ref{iota:morita} that $\sign^\M_\iota$ is uniquely defined up to a choice of sign.
We can arbitrarily choose the sign of the signature of a form at 
each ordering $P$. See for instance Remark~\ref{rem:scal} for a way to change sign using Morita equivalence (scaling).

A more intrinsic definition is therefore desirable, in particular when considering the total signature map of a hermitian form $X_F\to\Z$ since such arbitrary changes of sign would prevent it from being continuous. We are thus led to define a signature that is independent of the choice of $\iota$ and $\M$.

\begin{lemma}\label{lem:spiff} 
Let   $P\in X_F\sm\Nil[A,\s]$. Let $\iota_1 : F \rightarrow L_1$ and $\iota_2 : F \rightarrow L_2$ be two real
closures of $F$ at $P$ and let $\M_1$ and $\M_2$ be two different Morita equivalences as in \eqref{eq:morita}. 
Let $h_0\in W(A,\s)$ be such that 
$\sign_{\iota_1}^{\M_1} h_0 \not=0$ and let $\delta_k \in \{-1,1\}$ be the sign of $\sign_{\iota_k}^{\M_k} h_0$ for $k=1,2$. 
Then
\[\delta_1 \sign_{\iota_1}^{\M_1} h = \delta_2 \sign_{\iota_2}^{\M_2} h,\]
for all $h\in W(A,\s)$.
\end{lemma}

\begin{proof} Let $\delta\in\{-1,1\}$ be as in Corollary~\ref{iota:morita}. Then, we have 
for all $h\in W(A,\s)$ that
$\sign_{\iota_2}^{\M_2} h = \delta \sign_{\iota_1}^{\M_1} h$
and in particular that 
$\sign_{\iota_2}^{\M_2} h_0 = \delta \sign_{\iota_1}^{\M_1} h_0$.
It follows that $\delta_1= \delta \delta_2$. Thus
$\delta_1 \sign_{\iota_1}^{\M_1} h = \delta \delta_2 \sign_{\iota_1}^{\M_1} h
 = \delta_2 \sign_{\iota_2}^{\M_2} h$.\qed
\end{proof}

We will show in Theorem~\ref{thm:main} that there exists a finite tuple $H=(h_1,\ldots, h_s)$ of diagonal hermitian forms of dimension one over $(A,\s)$ such that for every $P\in X_F\sm\Nil[A,\s]$ there exists 
$h_0\in H$ such that $\sign_\iota^\nc h_0\not=0$.

\begin{definition}\label{def:s-sign} 
Let $h\in W(A,\s)$ and let $P\in X_F$. 
We define the \emph{$H$-signature of $h$ at $P$} as  follows:
If $P\in \Nil[A,\s]$, define
$\sign_P^H h:=0$.
If $P\not\in \Nil[A,\s]$,  let $i\in \{1,\ldots, s\}$ be the least integer such that $\sign_\iota^\M h_i\not=0$ (for any $\iota$ and $\M$, cf. Corollary~\ref{iota:morita}), let $\delta_{\iota,\M} \in \{-1,1\}$ be the sign of $\sign_\iota^\M h_i$ and define
\[\sign_P^H h:= \delta_{\iota,\M} \sign_\iota^\M h.  \]
\end{definition}

Lemma~\ref{lem:spiff} shows that this definition is independent of the choice of $\iota$ and $\M$
(but it does depend on the choice of $H$). 

A choice of Morita equivalence which is convenient for  computations of signatures is given by \eqref{comp:morita} with $(A,\s)$ replaced by $(A\ox_F F_P, \s\ox\id)$. We denote this Morita equivalence by $\CN$ and now describe the induced isomorphisms of Witt groups:
\begin{small}
\begin{equation}\label{seq5}
\begin{split}
\xymatrix@R=1ex{
W(A\ox_F F_P, \s\ox \id)\ar[r]^--{\xi_P^*}& W(M_m(D_P), \ad_{\vf_P}) \ar[r]^--{\text{scaling}} &W(M_m(D_P), {\vt_P}^t)
\ar[r]^--{\text{collapsing}} & W(D_P, \vt_P)\\
h\ox F_P\ar@{|->}[r] & \xi_P^*(h\ox F_P)\ar@{|->}[r] & \Phi_P^{-1} \xi_P^*(h\ox F_P) \ar@{|->}[r] & \CN(h\ox F_P),
}
\end{split}
\end{equation}
\end{small}
where $h$ is a hermitian form over $(A,\s)$, $P\in X_F\sm \Nil[A,\s]$ (so that $\ve_P=1$),
$\xi_P^*$ is the group isomorphism  induced by some fixed isomorphism
\[\xi_P: (A\ox_F F_P, \s\ox\id)\simtoo (M_m(D_P), \ad_{\vf_P}),\]
and $\Phi_P$ is the Gram matrix of the form $\vf_P$. Observe that $\sign \vf_P$ can be computed as in Section~\ref{sec:special}.

\begin{lemma}\label{lem5.5} Let $P \in X_F\sm \Nil[A,\s]$, let $\iota:F \to F_P$ be a real closure of $F$ and let $\CN$ and $\vf_P$ be as above. Then $\sign^\CN_\iota \<1\>_\s =\sign \vf_P$.
\end{lemma}

\begin{proof} 
We extend scalars from $F$ to $F_P$ via $\iota$, $\<1\>_\s \mapstoo \<1\>_\s\ox F_P=
 \<1\ox 1\>_{\s\ox\id}$ and push $ \<1\ox 1\>_{\s\ox\id}$  through the sequence $\eqref{seq5}$, 
\[\<1\ox 1\>_{\s\ox\id} \mapstoo \xi^*_P (\<1\ox 1\>_{\s\ox\id}) = \<\xi_P(1\ox 1)\>_{\ad_{\vf_P}}
=\<I_m\>_{\ad_{\vf_P}}\mapstoo \Phi_P^{-1}\<I_m\>_{\ad_{\vf_P}}=
\<\Phi_P^{-1}\>_{{\vt_P}^t}.\]
(Note that  $\xi_P (1\ox 1)=I_m$, the $m\x m$-identity matrix in $M_m(D_P)$ since  $\xi_P$ is an algebra homomorphism.) By collapsing, the matrix $\Phi_P^{-1}$ now corresponds to a quadratic form over $F_P$, a hermitian form over $(F_P(\sqrt{-1}),\bbar)$ or a hermitian form over $(H_P, \bbar)$. In either case $\Phi_P^{-1}$ is congruent to $\Phi_P$.  Thus $\sign^\CN_\iota \<1\>_\s =\sign \vf_P$.\qed
\end{proof}

\begin{remark}  It follows from Lemma~\ref{lem5.5} that the signature defined in \cite[\S3.3, \S3.4]{BP2}  is actually $\sign^{H}_P$ with $H=(\<1\>_\s)$.
It is now clear that this signature cannot be computed when $\sign_\iota^\CN \<1\>_\s=0$, i.e. when $\s\ox \id_{F_P}\simeq \ad_{\vf_P}$ is hyperbolic.  In contrast, if we take $H=(h_1,\ldots, h_s)$, as described before Definition~\ref{def:s-sign} we are able to compute
the signature in all cases. Note that we may choose $h_1=\<1\>_\s$, so that Definition~\ref{def:s-sign}  generalizes the definition of signature in \cite[\S3.3, \S3.4]{BP2}.
\end{remark}

\begin{example} Let $(A,\s) = (M_4(\R), \ad_\vf)$, where $\vf=\<1,-1,1,-1\>$. Then $\sign \vf=0$. Consider the dimension one hermitian forms
\[ h= \left\< \left(\begin{smallmatrix}
1& & & \\
 & 1 & & \\
& & -1 & \\
& & & 1
\end{smallmatrix}\right)\right\>_\s,\  
h_1= \left\< \left(\begin{smallmatrix}
1& & & \\
 & 1 & & \\
& & 1 & \\
& & & 1
\end{smallmatrix}\right)\right\>_\s, \text{ and }
h_2= \left\< \left(\begin{smallmatrix}
1& & & \\
 & -1 & & \\
& & 1 & \\
& & & -1
\end{smallmatrix}\right)\right\>_\s
\]
over $(A,\s)$. Then $\sign^\CN h=-2$, $\sign^\CN h_1=0$ and $\sign^\CN h_2=4$, where we suppressed the index $\iota$ since $\R$ is real closed. Let $H_1=(h_1)$ and $H=(h_1,h_2)$, then 
$\sign^{H_1}h$ is not defined, whereas $\sign^H h=  -2$. Observe that taking $H=(h_1, -h_2)$ instead would result in $\sign^H h=  2$.
\end{example}

\begin{remark}\label{rem:scal} 
Let $a\in A^\x$ be such that  $\s(a)=\ve a$ with $\ve\in\{-1,1\}$.
The Morita equivalence \emph{scaling by $a$},
$\Herm(A,\s) \too \Herm_\ve(A, \Int(a)\circ\s) ,\ h \mapstoo ah$
induces an isomorphism 
$\zeta_a: W(A,\s) \too W_\ve(A, \Int(a)\circ\s) ,\ h \mapstoo ah$.
It is clear that
\[\sign_\iota^\M h = \sign_\iota^{\M\circ (\zeta_a^{-1}\ox \id)} ah.\]

Consider the special case where $a\in F^\x$. Thus $\ve=1$ and $\Int(a)\circ \s =\s$. Assume that
$a<_P 0$. Then 
$\sign_\iota^\M \zeta_a(h) = \sign \M(ah\ox F_P)= \sign \M(-h\ox F_P) = -\sign_\iota^\M h$, 
where the last equality follows from Proposition~\ref{rem:morsign}$(iii)$.
The same computation shows that $\sign^H_P \zeta_a(h) = - \sign^H_P h$ for any choice of $H$.
Thus, scaling by $a$ changes the sign of the signature, which is contrary to what is claimed in \cite[p. 662]{BP2}.
\end{remark}

\begin{remark}\label{rem:Mchoice} 
For any choice of $H$, $P$ and $\iota$ as in Definition~\ref{def:s-sign}, there exists a Morita equivalence $ \M'$ such that $\sign_P^H h=\sign_\iota^{\M'} h$ for any $h\in W(A,\s)$ (i.e. such that $\sign_\iota^{\M'} h_i>0$ with $h_i$ as in Definition~\ref{def:s-sign}). Indeed, for $\M$ as in Definition~\ref{def:s-sign}, it suffices to take $\M'=\delta_{\iota,\M} \M$.
\end{remark}

It remains to be shown that a tuple $H$ as described just before Definition~\ref{def:s-sign} 
does exist. In order to reach this conclusion we first need to develop more theory in Sections~\ref{sec:inv} and \ref{ktf-M}.

\section{Signatures of Involutions}\label{sec:inv}
Let $(A,\s)$ be an $F$-algebra with involution of any kind with centre a field $K$. 
Consider the \emph{involution trace  form}
$T_\s: A\x A \too K, (x,y)\mapstoo \Trd_A(\s(x)y)$,
where $\Trd_A$ denotes the reduced trace of $A$. If $\s$ is of the first kind, $T_\s$ is a symmetric bilinear form over $F$. 
If $\s$ is of the second kind, $T_\s$ is a hermitian form over $(K,\s|_{K})$. Let $P\in X_F$. The 
\emph{signature of the involution $\s$ at $P$} is defined by
$\sign_P \s:= \sqrt{\sign_P T_\s}$
and is a nonnegative integer, since $\sign_P T_\s$ is always a square; cf.
Lewis and Tignol \cite{LT} for involutions of the first kind and Qu\'eguiner \cite{Q} for involutions of the second kind.
We call the involution $\s$ \emph{positive at $P$} if $\sign_P \s\not=0$.

\begin{example} \mbox{} \label{trace}
\begin{enumerate}[$(i)$]

\item Let $(A,\s)=(M_n(F), t)$, where $t$ denotes transposition. Then $T_\s\simeq n^2\x \<1\>$. Hence  $\sign_P\s=n$ for all $P \in X_F$.

\item Let $(A,\s)=((a,b)_F, \bbar)$, where $\bbar$ denotes quaternion conjugation. Then $T_\s \simeq \<2\>\ox \<1,-a,-b,ab\>$.  Hence
$\sign_P \s=2$ for all $P\in X_F$ such that $a<_P 0, b<_P 0$ and $\sign_P=0$ for all other $P\in X_F$. 

\item Let $(A,\s)=(F(\sqrt{d}), \bbar)$, where $\bbar$ denotes conjugation. 
Then $T_\s\simeq  \<1\>_\s$. We have $\sign_P \<1\>_\s = \frac{1}{2} \sign_P \<1,-d\>$, cf.
\cite[Chapter 10, Examples 1.6(iii)]{Sch}.
Hence  $\sign_P\s=0$ for all $P \in X_F$ such that $d>_P0$ and
$\sign_P\s=1$ for all $P \in X_F$ such that $d<_P0$. 

\end{enumerate}
\end{example}

\begin{remark}\label{rem3.8}
 Let $(A,\s)$ and $(B,\tau)$ be two $F$-algebras with  involution. 
\begin{enumerate}[$(i)$]
\item Consider the tensor product $(A\ox_F B, \s\ox\tau)$.
Then $T_{\s\ox\tau}\simeq T_\s \ox T_\tau$ and so $\sign_P(\s\ox\tau)=(\sign_P \s)(\sign_P \tau)$ for all $P \in X_F$.

\item If $(A,\s) \simeq (B,\tau)$, then $T_\s \simeq T_\tau$ so that $\sign_P\s= \sign_P \tau$ for all $P\in X_F$.
\end{enumerate}

\end{remark}

\begin{remark} \label{rem3.9}
Pfister's local-global principle holds for algebras with involution $(A,\s)$ and also for hermitian forms $h$ over such algebras,  \cite{LU1}.
\end{remark}

\begin{remark}\label{rem4.6} The map $\sign \s$ is continuous from $X_F$
(equipped with the Harrison topology, see \cite[Chapter VIII 6]{Lam} for a
definition) to $\Z$ (equipped with the discrete topology). Indeed: define the map
$\sqrt{\phantom{x}}$ on $\Z$ by setting $\sqrt{k}=-1$ if $k$ is not a
square in $\Z$. Since $\Z$ is equipped with the discrete topology, this map
is continuous. Since $T_\s$ is a symmetric bilinear form or a hermitian form over $(K,\s|_K)$, 
the map $\sign T_\s$ is
continuous from $X_F$ to $\Z$ (by \cite[VIII, Proposition 6.6]{Lam} and \cite[Chapter~10, Example~1.6(iii)]{Sch}). Thus, by
composition, $\sign \s = \sqrt{\sign T_\s}$ is continuous from $X_F$ to
$\Z$.
\end{remark}

\begin{lemma}\label{lem5.6} 
Let $P\in X_F$. If $P\in \Nil [A,\s]$, then 
$\sign_P\s = \sign \vf_P=0$. 
Otherwise,
$\sign_P \s = \lambda_P\,  |\!\sign \vf_P|$,
where $\lambda_P=1$ if $(D_P,\vt_P) = (F_P,\id_{F_P})$ or $(D_P,\vt_P) = (F_P(\sqrt{-1}), -)$  and $\lambda_P=2$ if $(D_P,\vt_P) = (H_P,\bbar)$.
\end{lemma}

\begin{proof}   This is a reformulation of \cite[Theorem~1]{LT} and part of its proof for involutions of the first kind and \cite[Proposition~3]{Q} for involutions of the second kind.\qed
\end{proof}

\begin{lemma}\label{h-adh} Let $(M,h)$ be a hermitian space over $(A,\s)$, let $P\in X_F$, let $\iota : F\to F_P$ be a real closure of $F$ at $P$ and let 
$\M$ be a Morita equivalence as in  \eqref{eq:morita}. 
 If $P\in \Nil [A,\s]$, then 
 $\sign_P\ad_h = \sign^\nc_\iota h=0$. 
 Otherwise,
$\sign_P \ad_h = \lambda_P\, |\! \sign_\iota^\nc h|$,
with $\lambda_P$ as defined in Lemma~\ref{lem5.6}. In particular,
$\sign_\iota^\nc h= 0$ if and only if $\sign_P \ad_h=0$.
\end{lemma}

\begin{proof} Assume first that $P\in \Nil[A,\s]$. Then $\sign_\iota^\nc h=0$. 
Consider the adjoint involution $\ad_h$ on $\End_A(M)$. Since $h$ is hermitian, $\s$ and $\ad_h$ are of the same type. Furthermore, $A$ and $\End_A(M)$ are Brauer equivalent by \cite[1.10]{BOI}. Thus $\Nil[A,\s] = \Nil[\End_A(M), \ad_h]$. By Lemma~\ref{lem5.6} we conclude that $\sign_P \ad_h=0$.

Next, assume that $P\in X_F\sm \Nil[A,\s]$. Without loss of generality we may replace $F$ by  $F_P$.  Consider a Morita equivalence
$\M':\Herm(A,\s) \too \Herm (D,\vt)$
with $(D,\vt) = (F,\id)$, $(D,\vt) = (H, \bbar)$ or $(D,\vt) = (F(\sqrt{-1}),-)$.  
Let $(N,b)$ be the hermitian space over $(D,\vt)$ corresponding to $(M,h)$ under $\M'$.
Then $\sign^{\M'} h=\sign b$. By \cite[Remark 1.4.2]{BP1} we have $(\End_A(M), \ad_h) \simeq (\End_D(N), \ad_b)$ so that $\sign \ad_h=\sign \ad_b$.
By  \cite[Theorem~1]{LT} and \cite[Proposition~3]{Q} we have $\sign \ad_b=\lambda\, |\!\sign b|$ with $\lambda=1$ if 
$(D,\vt)= (F,\id)$ or $(D,\vt)= (F(\sqrt{-1}),-)$ and $\lambda=2$ if $(D,\vt) = (H, \bbar)$. We conclude that $\sign \ad_h=\lambda\, |\!\sign^{\M'} h|=\lambda\, |\!\sign^{\M} h|$, where the last equality follows from Corollary~\ref{iota:morita}.\qed
\end{proof}

\section{The Knebusch Trace Formula for $M$-Signatures}\label{ktf-M}

We start with two preliminary sections in order not to overload the proof of Theorem~\ref{ktf} below.

\subsection{Hermitian Forms over a Product of Rings with
Involution}\label{product-rings}
Let 
$(A,\s) = (A_1,\s_1) \times \cdots \times (A_t, \s_t)$,
where $A,A_1,
\ldots, A_t$ are rings and $\s, \s_1, \ldots, \s_t$ are involutions. We write an
element $a \in A$ indiscriminately as $(a_1, \ldots, a_t)$ or $a_1 + \cdots +
a_t$ with $a_i \in A_i$ for $i=1, \ldots, t$. Writing $1_A = (e_1, \ldots,
e_t)$, the elements $e_1, \ldots, e_t$ are central idempotents, and 
the coordinates of $a\in A$ are given by
\[A \too A_1 \x \cdots \x A_t,\ a \mapstoo (ae_1, \ldots, ae_t).\] 
Note that $e_ie_j=0$ whenever $i\not=j$. We assume that
$\s(1)=1$ and thus $\s(e_i)=e_i$ for $i=1, \ldots, t$.

Let $M$ be an $A$-module and let $h : M \x M \rightarrow A$ be a hermitian form
over $(A,\s)$. Following \cite[Proof of Lemma~1.9]{KRW}
we can write 
\[M  = \coprod_{i=1}^t Me_i,\ m = (me_1, \ldots, me_t),\]
where $\coprod_{i=1}^t Me_i$ is the $A$-module
with set of elements $\prod_{i=1}^t Me_i$, whose sum is defined coordinate by
coordinate and whose product is defined by $(m_1e_1, \ldots, m_te_t) a = (m_1e_1a_1, \ldots, m_te_ta_t)$ for $m_1, \ldots, m_t \in M$ and $a \in A$.

Define $h_i = h|_{Me_i}$. Then $h_i(xe_i,ye_i) = \sigma(e_i)h(x,y)e_i =
h(x,y)e_i^2 = h(x,y)e_i$ and $h_i
: Me_i \x Me_i \rightarrow A_i$ is a hermitian form over $(A_i, \s_i)$. We also
have
\begin{equation*}
\begin{split}
  h(xe_1 + \cdots + xe_t, ye_1 + \cdots + ye_t) &= \sum_{i,j=1}^t
      h(xe_i,ye_j)e_ie_j \\
    &= \sum_{i=1}^t h(xe_i,ye_i)e_i
\end{split}
\end{equation*}
which proves that $h = h_1 \perp \ldots \perp h_t$.

\subsection{Algebraic Extensions and Real Closures}\label{alg-ext}
We essentially follow \cite[Chapter~3, Lemma~2.6, Lemma~2.7, Theorem~4.4]{Sch}.

Let $P\in X_F$ and let $F_P$ denote a real closure of $F$ at $P$. 
Let $L$ be a finite extension of $F$.
Writing $L =
F[X]/(R)$ for some $R \in F[X]$ and $R = R_1 \cdots R_t$ as a product of pairwise
distinct irreducibles in $F_P[X]$ with $\deg R_1 = \cdots = \deg R_r = 1$ and
$\deg R_{r+1} = \cdots = \deg R_t = 2$, we obtain canonical $F_P$-isomorphisms
$L\ox_F F_P \simeq F_P[X]/(R_1 \cdots R_t)$
and
\begin{equation}\label{eq:LFP}
  L \ox_F F_P  \stackrel{\omega}{\too} F_1 \x \cdots \x F_t,
\end{equation}
where $F_i = F_P[X]/(R_i)$ is a real closed field for $1 \le
i \le r$ and is algebraically closed for $r+1 \le i \le t$. We write $1 = (e_1,
\ldots, e_t)$ in $F_1 \x \cdots \x F_t$ and define $ \omega_i(x) = \omega(x)e_i$
for $x \in L \ox_F F_P$, the projection of $\omega(x)$ on its $i$-th coordinate.

Let $\iota_0 : F_P \rightarrow L \ox_F F_P$ be the canonical inclusion. Then
$ \omega_i \circ \iota_0 : F_P \rightarrow F_i$ is an isomorphism of fields and
of $F_P$-modules for $1 \le i \le r$. In particular, for $1 \le i \le r$, $F_i$
is naturally an $F_P$-module of dimension one, and it is easily seen that
$\Tr_{F_i/F_P} = ( \omega_i \circ \iota_0)^{-1}$ and so $\Tr_{F_i/F_P}$
is an isomorphism of fields.

Let $\iota_1 : L \rightarrow L \ox_F F_P$ be the canonical inclusion. Then
$ \omega_i \circ \iota_1 : L\to F_i$ ($i=1, \ldots, r$) denote the $r$
different embeddings of ordered fields corresponding to the orderings $Q_i$ on
$L$ that extend $P$.
In other words, if $\{Q_1, \ldots, Q_r\}$ are the different extensions of $P$ to
$L$, then for every $1 \le i \le r$, the map 
\[\xymatrix{
  L \ar[r]^--{\iota_1} & L \ox_F F_P \ar[r]^--{\omega_i} & F_i
  }\]
is a real closure of $L$ at $Q_i$. Since $\Tr_{F_i/F_P}$ is an isomorphism of
fields, it follows that the map
\[\xymatrix{
  L \ar[r]^--{\iota_1} & L \ox_F F_P \ar[r]^--{\omega_i} & F_i \ar[r]^--{\Tr_{F_i/F_P}} &
  F_P
  }\]
is also a real closure of $L$ at $Q_i$.

\subsection{The  Knebusch Trace Formula} 
Let $(A,\s)$ be an $F$-algebra with involution. 
Let $L/F$ be a finite extension. The trace  $\Tr_{L/F}: L\to F$ induces an $A$-linear homomorphism
\[\Tr_{A\ox_F L}  =\id_A \ox \Tr_{L/F} : A\ox_F L \too A\]
which induces a group homomorphism (transfer map)
\[\trs_{A\ox_F L}: W(A\ox_F L, \s\ox\id) \too W(A,\s),\ (M,h)\mapstoo (M, \Tr_{A\ox_F L} \circ h),\] 
cf. \cite[p. 362]{B-L}.

The following theorem is an extension of a result due to Knebusch \cite[Proposition~5.2]{K}, \cite[Chapter~3, Theorem~4.5]{Sch} to $F$-algebras with involution. The proof follows the general lines of Knebusch's original proof.

\begin{theorem}\label{ktf} Let $P\in X_F$. Let $L/F$ be a finite extension of
ordered fields and let $h$ be a hermitian form over $(A\ox_F L, \s\ox\id)$. Fix
a real closure $\iota : F \rightarrow F_P$ and a Morita equivalence $\nc$ as in
\eqref{eq:morita}. Then, with notation as in Section~\ref{alg-ext},
\[\sign_{\iota}^\nc (\trs_{A\ox_F L} h) =\sum_{i=1}^r
\sign_{\omega_i \circ \iota_1}^{(\omega_i \circ \iota_0)(\nc)} h.\]
\end{theorem}
\begin{proof} 
By definition of signature we have
\begin{equation}\label{sign_1}
\sign_\iota^\nc (\trs_{A\ox_F L} h) = \sign \nc [ (\trs_{A\ox_F L} h)\ox_F
F_P].
\end{equation}
Consider the commutative diagram
\begin{equation*}
\xymatrix{
L \ar[rr]^--{\Tr_{L/F} } \ar[d]^--{\ox_F F_P}& & F\ar[d]^--{\ox_F F_P}\\
L\ox_F F_P \ar[rr]^--{\Tr_{L\ox F_P/F_P}} &  &F_P\\
}
\end{equation*}
It induces a commutative diagram
\begin{equation*}
\xymatrix{
A\ox_F L \ar[rr]^--{\Tr_{A\ox L} } \ar[d]^--{\ox_F F_P}& & A\ar[d]^--{\ox_F F_P}\\
A\ox_F L\ox_F F_P \ar[rr]^--{\id_A\ox \Tr_{L\ox F_P/F_P}} &  &A\ox_F F_P\\
}
\end{equation*}
which in turn induces a commutative diagram of Witt groups
\begin{equation*}
\xymatrix{
W(A\ox_F L,\s\ox\id) \ar[rr]^--{\trs_{A\ox L} } \ar[d]& & \ar[d] W(A,\s) \\
W(A\ox_F L\ox_F F_P, \s\ox\id\ox\id) \ar[rr]^--{(\id_A\ox \Tr_{L\ox F_P/F_P})^* } 
&  &W(A\ox_F F_P, \s\ox\id)\\
}
\end{equation*}
where the vertical arrows are the canonical restriction maps. Thus
\begin{equation}\label{sign_2}
\sign {\nc} [ (\trs_{A\ox_F L} h)\ox_F F_P]= 
\sign {\nc} [ (\id_A\ox \Tr_{L\ox F_P/F_P})^* (h\ox_F F_P)].
\end{equation}
With reference to Section~\ref{alg-ext}
consider the diagram 
\begin{equation*}
\xymatrix{
A\ox_F L \ox_F F_P \ar[r]^--{\id_A \ox \omega} \ar[d]^--{\id_A  \ox_F  \Tr_{L\ox F_P /F_P}}&  A  \ox_F (F_1\x\cdots \x F_t) \ar[r]^--{\sim}     \ar[d]^--{\id_A \ox  \Tr_{F_1\x\cdots\x F_t /F_P}  } & (A\ox_F F_1)\x \cdots \x (A\ox_F F_t)
\ar[d]^--{\sum_{i=1}^t \id_A \ox \Tr_{F_i/F_P} } \\
A\ox_F  F_P \ar[r]^--{\id} &A\ox_F  F_P \ar[r]^--{\id}  &A\ox_F F_P\\
}
\end{equation*}
where commutativity of the first square follows from the isomorphism
\eqref{eq:LFP} of $F_P$-algebras, whereas commutativity of the second square
follows from \cite[p.~137]{B-alg8}. We push $h\ox F_P$ through the induced
commutative diagram of Witt groups:
\begin{equation*}
\xymatrix{
h\ox F_P \ar@{|->}[r]^--{(\id_A\ox\omega)^*} \ar@{|->}[d] & h' \ar@{|->}[r] & h'_1 \perp \ldots \perp h'_t \ar@{|->}[d]\\
(\id_A\ox \Tr_{L\ox F_P / F_P})^* (h\ox F_P) \ar@{|->}[rr]^--{\id} & & \sum_{i=1}^t (\id_A\ox \Tr_{F_i / F_P})^*(h'_i)
}
\end{equation*}
where the image of $h'$ equals the orthogonal sum $h'_1\perp \ldots \perp h'_t$ by Section~\ref{product-rings}.
Thus
\begin{equation}\label{sign_3}
 \sign {\nc} [ (\id_A\ox \Tr_{L\ox F_P/F_P})^* (h\ox_F F_P)] = \sum_{i=1}^t
 \sign {\nc}[(\id_A\ox \Tr_{F_i / F_P})^*(h'_i)].
\end{equation}
We have to consider the following two cases:

Case 1: Assume that $1\leq i \leq r$. Observe that $h'_i = (\id_A \ox( \omega_i \circ
\iota_1))^*(h)$. Then
\begin{equation*}
  \sign {\nc}[(\id_A\ox \Tr_{F_i / F_P})^*(h'_i)] 
    = \sign {\nc}[(\id_A\ox (\Tr_{F_i / F_P}))^*\circ (\id_A \ox (\omega_i \circ
      \iota_1))^*(h)].
\end{equation*}
The form $(\id_A \ox (\omega_i \circ \iota_1))^*(h)$ is defined over $A \ox_F
F_i$, and the commutative diagram
\[\xymatrix@R=1ex{
&  F_P\ar[dd]^{\omega_i \circ \iota_0}\\
  F_i \ar[ru]^{\Tr_{F_i / F_P}} \ar[dr]_{\id} &  \\
&  F_i\\
  }\]
together with Proposition~\ref{prop:iota} gives
\begin{multline}\label{sign_5}
  \sign {\nc}[(\id_A\ox (\Tr_{F_i / F_P}))^*\circ  (\id_A \ox (\omega_i \circ
      \iota_1))^*(h)]  \\
      =\sign (\omega_i \circ \iota_0)(\M)[(\id_A \ox (\omega_i \circ
      \iota_1))^*(h)] 
\end{multline}

Case 2: Assume that $r+1\leq i \leq t$.  Since $F_i$ is algebraically closed, it follows from
Morita theory that the Witt group $W(A\ox_F F_i, \s\ox\id)$ is torsion
 and so $h'_i$ is a
torsion form. Therefore $(\id_A\ox \Tr_{F_i / F_P})^*(h'_i)$ is also a torsion
form and thus has signature zero. 

We conclude that equations \eqref{sign_1}-\eqref{sign_5} yield the Knebusch trace formula.\qed  
\end{proof}

\section{Existence of Forms with Nonzero Signature}\label{choose_sign}

\begin{theorem}\label{thm2.4}
 Let $(A,\s)$ be an $F$-algebra with involution and let $P\in X_F\sm \Nil[A,\s]$. 
 Fix
a real closure $\iota : F \rightarrow F_P$ and a Morita equivalence $\nc$ as in
\eqref{eq:morita}.
 There exists a hermitian form $h$ over $(A,\s)$ such that $\sign_\iota^\nc h\not=0$.
\end{theorem}

\begin{proof} Let $P\in X_F\sm \Nil[A,\s]$. Then $A\ox_F F_P \simeq  M_m(D_P)$ for some $m\in \N$, where
$D_P=F_P$, $H_P$ or $F_P(\sqrt{-1})$ if $\s$ is orthogonal, symplectic or
unitary, respectively. In each case there exists a positive involution $\tau$ on $M_m(D_P)$
(namely, transposition, conjugate transposition and quaternion conjugate transposition, respectively, cf. Example~\ref{trace}).
By Lemma~\ref{h-adh} the hermitian form $\<1\>_\tau$ over $(M_m(D_P), \tau)$ has nonzero signature since $\tau$ is the adjoint involution of the form $\<1\>_\tau$. After scaling we obtain a dimension one hermitian form $h_0$ over $(A\ox_F F_P,
\s\ox\id)$ such that $\sign \M(h_0)\not=0$
by Proposition~\ref{prop:morita}.  The form $h_0$ is already defined over a finite extension $L$ of $\iota(F)$, contained in $F_P$. Thus we can consider $h_0$ as a form over $(A\ox_F L, \s\ox\id)$ and we have  $\sign \M(h_0\ox F_P)\not=0$. In other words, if  $Q_1$ is the ordering 
$L\cap F_P^{\x 2}$ on $L$, then for any real closure $\kappa_1: L\to L_1$  of $L$ at $Q_1$ and any Morita equivalence $\M_1$ as in \eqref{eq:morita}, but starting from $\Herm (A\ox_F L_1, \s\ox\id)$, we have 
$\sign_{\kappa_1}^{\M_1} h_0\not=0$
by Corollary~\ref{iota:morita}. 

Let $X=\{Q\in X_L \mid
P\subset Q\}$. 
By \cite[Chapter~3, Lemma~2.7]{Sch}, $X$ is finite, say $X=\{Q_1,Q_2,\ldots, Q_r\}$. Thus
there exist $a_2,\ldots, a_r\in L^\x$ such that
\[\{Q_1\} =H(a_2,\ldots, a_r)\cap X.\]
Consider the Pfister form $q:=\<\!\<a_2,\ldots, a_r\>\!\>=\<1,a_2\> \ox \cdots \ox \<1,a_r\>$. Then $\sign_{Q_1} q=2^{r-1}$ and $\sign_{Q_\ell}q=0$ for $\ell\not=1$. It follows that 
$\sign^{\M_1}_{\kappa_1} (q\ox h_0)=\sign_{Q_1}q\cdot \sign^{\M_1}_{\kappa_1} h_0\not=0$ and $\sign^{\M_\ell}_{\kappa_\ell} 
(q\ox h_0)=0$ for $\ell\not=1$,
where $\kappa_\ell: L\to L_\ell$ is any real closure of $L$ at $Q_\ell$ and  $\M_\ell$ is any Morita equivalence as in \eqref{eq:morita}, but starting from $\Herm (A\ox_F L_\ell, \s\ox\id)$.

Now $\trs_{A \ox_F L}(q\ox h_0)$ is a hermitian form over $(A,\s)$. By the trace formula, Theorem~\ref{ktf}, we have
\[\sign_{\iota}^\nc \bigl(\trs_{A\ox_F L} (q\ox h_0)\bigr) =\sum_{i=1}^r
\sign_{\omega_i \circ \iota_1}^{(\omega_i \circ \iota_0)(\nc)} (q\ox h_0)=\sign_{\omega_1 \circ \iota_1}^{(\omega_1 \circ \iota_0)(\nc)} (q\ox h_0)\not=0.\]
Taking $h:= \trs_{A\ox_F L} (q\ox h_0)$ proves the theorem.\qed
\end{proof}

\begin{corollary} \label{cor4.2}
Let $(A_1,\s_1)$ and $(A_2,\s_2)$ be $F$-algebras with involution of the same
type.  Assume that $A_1$ and $A_2$ are Brauer equivalent. Let $P\in X_F$, let $\iota: F \to F_P$ be a real closure of $F$ at $P$ and let
$\M_\ell$ be any Morita equivalence as in \eqref{eq:morita}, but starting from $\Herm (A_\ell\ox_F F_P, \s_\ell\ox\id)$ for $\ell=1,2$.
Then the following statements are equivalent:
\begin{enumerate}[$(i)$]
  \item $\sign_\iota^{\M_1} h=0$ for all hermitian forms $h$ over $(A_1,\s_1)$;
  \item $\sign_\iota^{\M_2} h=0$ for all hermitian forms $h$ over $(A_2,\s_2)$;
  \item $\sign_P \vt=0$ for all involutions $\vt$ on $A_1$ of the
    same type as $\s_1$;
  \item $\sign_P \vt=0$ for all involutions $\vt$ on $A_2$ of
    the same type as $\s_2$.
\end{enumerate}
\end{corollary} 

\begin{proof}  
By Theorem~\ref{thm2.4}, the first two statements are equivalent  to $P\in \Nil[A_1,\s_1]=\Nil[A_2,\s_2]$. Thus $(i)\iff (ii)$.

$(i)\implies (iii)$ Let $\vt$ be as in $(iii)$. Then $\vt=\ad_{\<1\>_\vt}$ and $\vt = \Int(a)\circ \s_1$ for some invertible $a\in \Sym(A_1,\s_1)$. Thus, with notation as in Remark~\ref{rem:scal} and using Lemma~\ref{h-adh} and Proposition~\ref{prop:morita} we have for any Morita equivalence $\M$,
\begin{align*}
\sign_P \vt&= \lambda_P | \sign_\iota^\M \<1\>_\vt| \\
&=\lambda_P | \sign_\iota^{\M\circ (\zeta_a\ox \id)} \zeta_a^{-1} (\<1\>_\vt)| \\
&= \lambda_P | \sign_\iota^{\M_1} \zeta_a^{-1} (\<1\>_\vt)| \\ 
&= \lambda_P | \sign_\iota^{\M_1} \<a^{-1}\>_{\s_1}| 
\end{align*}
which is zero by the assumption. 

$(ii)\implies (iv)$: This is the same proof as $(i)\implies (iii)$ after replacing $(A_1,\s_1)$ by $(A_2,\s_2)$.

For the remainder of the proof we may assume
without loss of generality  that $A_2$ is a division algebra and that $A_1\simeq M_m(A_2)$ for some $m\in \N$.

$(iii)\implies (iv)$: Let $\vt$ be any involution on $A_2$. By the assumption, $\sign_P (t\ox \vt)=0$, where $t$ denotes the transpose involution. Since $t$ is a  positive involution, it follows that $\sign_P \vt=0$, cf. Remark~\ref{rem3.8}.

$(iv) \implies (ii)$: The assumption implies that $\sign_\iota^{\M_2} h=0$ for
every hermitian form $h$ of dimension $1$ over $(A_2,\s_2)$, which implies $(ii)$
since all hermitian forms over $(A_2,\s_2)$ can be diagonalized.\qed 
\end{proof}

\begin{remark}\label{rem:signon} 
Note that statement $(iii)$ in Corollary~\ref{cor4.2} is equivalent to:  $\sign_\iota^{\M_1} h=0$ for all diagonal hermitian forms $h$ of dimension one over $(A_1,\s_1)$. \end{remark}

\begin{theorem}\label{thm:main} 
Let $(A,\s)$ be an $F$-algebra with involution. There exists a finite set $H=\{h_1,\ldots, h_s\}$ of diagonal hermitian forms of dimension one over $(A,\s)$ such that for every $P\in X_F\sm\Nil[A,\s]$, real closure $\iota : F \rightarrow F_P$ and  Morita equivalence $\M$ as in
\eqref{eq:morita}
there exists $h\in H$ such that $\sign_\iota^{\M} h\not=0$.
\end{theorem}

\begin{proof}  For every $P\in X_F$, choose a real closure $\iota_P : F \rightarrow F_P$ and  a Morita equivalence $\M_P$ as in
\eqref{eq:morita}. By Corollary~\ref{iota:morita} we may assume without loss of generality and for the sake of simplicity
that the map $\iota_P$ is an inclusion.   
The algebra $A\ox_F F_P$ is isomorphic to a matrix
algebra over $D_P$, where $D_P\in \{F_P, H_P, F_P(\sqrt{-1}), F_P\x F_P\}$.
There is a finite extension $L_P$ of $F$, $L_P\subset F_P$ such that $A\ox_F L_P$
is isomorphic to a matrix algebra over $E_P$, where $E_P\in \{L_P, (-1,-1)_{L_P},
L_P(\sqrt{-1}), L_P\x L_P\}$, and $P$ extends to $L_P$. Let
\[U_P:=\{ Q\in X_F \mid Q \text{ extends to } L_P\}.\]
Since $P\in U_P$ we can write
$X_F = \bigcup_{P\in X_F} U_P$.
We know from \cite[Chapter~3, Lemma~2.7, Theorem~2.8]{Sch} that $\sign_{Q} (\trs_{L_P/F} \<1\>)$ is the number of extensions of $Q$ to $L_P$. Thus,
\[U_P = \Bigl( \sign (\trs_{L_P/F} \<1\>)\Bigr)^{-1} (\{1,\ldots,k\}),\]
where $k$ is the dimension of the quadratic form $\trs_{L_P/F} \<1\>$, and so
$U_P$ is clopen in $X_F$ (and in particular compact). Therefore, and since $X_F$ is compact, there exists a finite number of orderings $P_1,\ldots, P_\ell$ in $X_F$ such that 
$X_F =\bigcup_{i=1}^\ell U_{P_i}$.
Now let $P\in \{P_1,\ldots, P_\ell\}$ and let $L_P$ be as before. By
Theorem~\ref{thm2.4} we have that for every $Q\in U_P \setminus \Nil[A,\s]$
there exists a hermitian form $h_Q$ over $(A,\s)$ such that  $\sign_{\iota_Q}^{\M_Q} h_Q \not=0$. 
By Corollary~\ref{cor4.2} and Remark~\ref{rem:signon} we may assume that $h_Q$ is diagonal of 
dimension one. 
Consider the total signature map
$\mu_Q: X_F\too \Z,\ P \mapstoo \sign_{\iota_P}^{\M_P} h_Q$.
Then
\begin{equation}\label{eq2.4-2}
U_P\setminus \Nil[A,\s] = \bigcup_{Q\in U_P\setminus \Nil[A,\s]} \mu_Q^{-1} (\Z\setminus \{0\}).
\end{equation}
Consider the continuous map
$\lambda_P: X_{L_P}\too X_F, R\mapstoo R\cap F$.
We have $Q\in U_P\setminus \Nil[A,\s]$ if and only if some extension $Q'$ of $Q$ to $L_P$ is in $X_{L_P}\setminus \Nil[A\ox_F L_P, \s\ox\id]$ (this follows from the fact that the ordered fields $(F,Q)$ and $(L_P,Q')$ have a common real closure)
if and only if $Q\in 
\lambda\bigl( X_{L_P}\setminus \Nil[A\ox_F L_P, \s\ox\id]\bigr)$. 

We observe that
$X_{L_P}\setminus \Nil[A\ox_F L_P, \s\ox\id]$ is clopen and compact since $\Nil[A\ox_F L_P, \s\ox\id]$ is 
either $\varnothing$ or the whole of $X_{L_P}$, which follows from the fact that $A\ox_F L_P$ is a matrix algebra over one of $L_P, (-1,-1)_{L_P},
 L_P(\sqrt{-1}), L_P\x L_P$.

Hence,
\[U_P\setminus \Nil[A,\s] = \lambda\bigl( X_{L_P}\setminus \Nil[A\ox_F L_P, \s\ox\id]\bigr)\]
is compact and thus closed. Thus $U_P\cap \Nil[A,\s]$ is open in $U_P$. Using \eqref{eq2.4-2} we can write
\[U_P=\bigl(U_P \cap \Nil[A,\s] \bigr) \cup \bigcup_{Q\in U_P\setminus \Nil[A,\s]}
\mu_Q^{-1} (\Z\setminus \{0\}).\]
Now $\mu_Q^{-1} (\Z\setminus \{0\}) = (\sign \ad_{h_Q})^{-1} (\Z\setminus \{0\})$ by Lemma~\ref{h-adh}, which
 is open since $\sign \ad_{h_Q}$ is continuous by Remark~\ref{rem4.6}. 
Thus, since $U_P$ is compact, there exist $Q_1,\ldots, Q_t \in U_P\setminus \Nil[A,\s]$ such that
\[U_P=\bigl(U_P \cap \Nil[A,\s] \bigr) \cup \bigcup_{i=1}^t
\mu_{Q_i}^{-1} (\Z\setminus \{0\}).\]
In other words, for every $Q \in U_P\setminus \Nil[A,\s]$ one of $\sign^{\M_{Q}}_{\iota_{Q}} h_{Q_i}$  $(i=1,\ldots, t)$ is nonzero.
Now let $H_P=\{h_{Q_1},\ldots, h_{Q_t}\}$. Letting $H=\bigcup_{i=1}^\ell H_{P_i}$ finishes the proof.\qed
\end{proof}

\begin{corollary} Let $(A,\s)$ be an $F$-algebra with involution. The set $\Nil[A,\s]$ is clopen in $X_F$.
\end{corollary}

\begin{proof} By Theorem~\ref{thm:main} we have $\Nil[A,\s] = \bigcap_{i=1}^s \{P\in X_F \mid \sign^{\M_P}_{\iota_P} h_i=0\}$. The result follows from Lemma~\ref{h-adh} and Remark~\ref{rem4.6}.\qed
\end{proof}

At this stage we have established all results that are needed for the definition of the $H$-signature in Definion~\ref{def:s-sign}. In the final two sections we show that the total $H$-signature  of a hermitian form is continuous and we reformulate the Knebusch trace formula in terms of the $H$-signature.

\section{Continuity of the Total $H$-Signature Map of a Hermitian Form}\label{sec:cont}

Let $h$ be a hermitian form over $(A,\s)$. With reference to Definition~\ref{def:s-sign} 
we denote by $\sign^H h$ the total $H$-signature map of $h$:
\[X_F\too \Z,\ P\mapstoo \sign_P^H h.\]

\begin{lemma}\label{lembluefrog} 
Let $H=(h_1,\ldots, h_s)$ be as in Definition~\ref{def:s-sign}.
There is a finite partition of $X_F$ into clopens
\[X_F= \Nil[A,\s] \dotcup \bigcupdot_{i=1}^\ell Z_i,\]
such that for every $i\in\{1,\ldots, \ell\}$ one of the total $H$-signature maps $\sign^H h_1,\ldots, \sign^H h_s$ is constant non-zero on $Z_i$.
\end{lemma}

\begin{proof}
For $r=1,\ldots, s$, let
\[Y_r:=\{ P\in X_F \mid \sign_P^H h_i=0,\ i=1,\ldots, r\}.\]
By Lemma~\ref{h-adh} we have
\[Y_r= \bigcap_{i=1}^r \{P\in X_F \mid \sign_P \ad_{h_i}=0\}.\]
Thus  each $Y_r$ is clopen. 

We have
$Y_0:=X_F \supseteq Y_1 \supseteq \cdots \supseteq Y_{s-1} \supseteq Y_s= \Nil[A,\s]$
and therefore,
\[X_F=(Y_0\sm Y_1) \dotcup (Y_1\sm Y_2) \dotcup \cdots \dotcup (Y_{s-1}\sm Y_s) \dotcup \Nil[A,\s].\]
Let $r\in \{0,\ldots,s-1\}$ and consider $Y_r\sm Y_{r+1}$. By the definition of $Y_1,\ldots, Y_s$
the map $\sign^H h_{r+1}$ is never $0$ on $Y_r\sm Y_{r+1}$. Furthermore,   $\sign^H h_{r+1}$  only takes a finite number of values $k_1,\ldots, k_m$. 

Now observe that 
there exists a $\lambda \in \{1,2\}$ such that 
\[\sign^H h_{r+1}=\frac{1}{\lambda} \sign \ad_{h_{r+1}}\]
on $Y_r\sm Y_{r+1}$. This follows from Lemma~\ref{h-adh} and Definition~\ref{def:s-sign} for $P\in Y_r\sm Y_{r+1}$.

Therefore,
\[ \bigl( \sign^H h_{r+1}\bigr)^{-1} (k_i) \cap (Y_r \sm Y_{r+1})=  \bigl(\sign  \ad_{h_{r+1}}\bigr)^{-1} (\lambda k_i)\cap (Y_r \sm Y_{r+1}),\]
which is clopen by Remark~\ref{rem4.6}. It follows that $Y_r \sm Y_{r+1}$ is covered by finitely many disjoint clopen sets on which the map 
$\sign^H h_{r+1}$ has constant non-zero value. The result follows since the sets $Y_r \sm Y_{r+1}$ for $r=0,\ldots, s-1$
form a partition of $X_F\sm \Nil [A,\s]$.\qed
\end{proof}

\begin{theorem} Let $h$ be a hermitian form over $(A,\s)$. The total $H$-signature of $h$, 
\[\sign^H h : X_F \too \Z,\ P\mapstoo \sign^H_P h\]
is continuous.
\end{theorem}

\begin{proof} We use the notation and the conclusion of Lemma~\ref{lembluefrog}. Since $\Nil[A,\s]$ and the sets $Z_i$ are clopen, it suffices to show that $(\sign^H h)|_{Z_i}$ is continuous for every $i=1,\ldots, \ell$. 

Let $i\in \{1,\ldots, \ell\}$,   $k_i\in \Z\sm\{0\}$ and  $j\in\{1,\ldots, s\}$
be such that $\sign^H h_j=k_i$ on $Z_i$. Let $k\in \Z$. Then
\begin{align*}
\bigl((\sign^H h)|_{Z_i}\bigr)^{-1} (k) &= \{P \in Z_i \mid \sign^H_P h = k\}\\
&= \{P \in Z_i \mid k_i \sign^H_P h = k_i k\}\\
&= \{P \in Z_i \mid k_i \sign^H_P h = k \sign^H_P h_j \}\\
&= \{P \in Z_i \mid  \sign^H_P  ( k_i \x h \perp -(k\x h_j))=0 \}.
\end{align*}
It follows from Lemma~\ref{h-adh} that
\[\bigl((\sign^H h)|_{Z_i}\bigr)^{-1} (k) =  \{P \in Z_i \mid  \sign_P \ad_{ k_i \x h \perp -(k\x h_j)}=0 \},\]
which is clopen by Remark~\ref{rem4.6}.\qed
\end{proof}

\section{The Knebusch Trace Formula for $H$-Signatures}\label{sec:ktf-H}

\begin{theorem} Let $H=(h_1,\ldots, h_s)$ be as in Definition~\ref{def:s-sign}. 
Let $P\in X_F$. Let $L/F$ be a finite extension of
  ordered fields and let $h$ be a hermitian form over $(A\ox_F L, \s\ox\id)$.
  Then, with $H\ox L:= (h_1\ox L,\ldots, h_s\ox L)$, we have 
  \[\sign^{H}_P (\trs_{A\ox_F L} h) =\sum_{P \subseteq Q \in X_L}
  \sign^{H \ox L}_Q h.\]
\end{theorem}

\begin{proof} 
  By Theorem \ref{ktf} (and using its notation), we know that 
  \[\sign_{\iota}^\nc (\trs_{A\ox_F L} h) =\sum_{i=1}^r \sign_{\omega_i \circ
  \iota_1}^{(\omega_i \circ \iota_0)(\nc)} h,\]
  for any $\iota$ and $\M$. Fix a real closure $\iota : F \rightarrow F_P$ and
  choose a Morita equivalence $\M$ such that $\sign_{\iota}^\M =
  \sign^{H}_P$ (cf. Remark~\ref{rem:Mchoice}). 
We only have to check that $\sign_{\omega_i \circ
  \iota_1}^{(\omega_i \circ \iota_0)(\nc)} = \sign^{H \ox L}_{Q_i}$ for 
  $i=1, \ldots, r$.

  By definition of $\M$, there is a $k \in \{1,
  \ldots, s\}$ such that $\sign_\iota^\M h_j = 0$ for $1 \le j \le k-1$
  and $\sign_\iota^\M h_k > 0$. To check that $\sign_{\omega_i \circ
  \iota_1}^{(\omega_i \circ \iota_0)(\nc)} = \sign^{H \ox L}_{Q_i}$ for
  $i=1, \ldots, r$, it suffices to check that $\sign_{\omega_i \circ
  \iota_1}^{(\omega_i \circ \iota_0)(\nc)}(h_j \ox L) = 0$ for $j=1,
  \ldots, k-1$ and $\sign_{\omega_i \circ \iota_1}^{(\omega_i \circ
  \iota_0)(\nc)}(h_k \ox L) > 0$. This follows from the fact that
  \begin{equation*}
  \sign_{\omega_i \circ \iota_1}^{(\omega_i \circ \iota_0)(\nc)}
  (h_\ell \ox L) = \sign_\iota^\M h_\ell  \text{ for every }1 \le \ell
  \le s,
  \end{equation*}
  which we verify in the remainder of the proof. 
  
  By definition,
  \[\sign_{\omega_i \circ \iota_1}^{(\omega_i \circ \iota_0)(\nc)}
  (h_\ell \ox L) = \sign (\omega_i \circ \iota_0)(\nc)[(\id_A \ox (\omega_i \circ
  \iota_1))^*(h_\ell \ox L)].\]
  Consider the commutative diagram
  \[\xymatrix@R=1ex{
          & F_P \ar[r]^--{\iota_0} & L \ox_F F_P \ar[r]^--{\omega_i} & F_i
          \ar[dd]^{\id} \\
    F \ar[ur]^\iota \ar[dr] & & & \\
          & L \ar[r]_--{\iota_1} & L \ox_F F_P \ar[r]_--{\omega_i} & F_i
    }\]
   Thus,  by Proposition \ref{prop:iota},
  \[ \sign (\omega_i \circ \iota_0)(\nc)[(\id_A \ox (\omega_i \circ \iota_1))^*
  (h_\ell \ox L)] = \sign (\omega_i \circ \iota_0)(\nc)[(\id_A \ox
  (\omega_i \circ \iota_0 \circ \iota))^*(h_\ell)].\]
  Finally the commutative diagram
   \[\xymatrix@R=1ex{
         & F_i \ar[dd]^{(\omega_i \circ \iota_0)^{-1}} \\
     F \ar[ur]^{\omega_i \circ \iota_0 \circ \iota} \ar[dr]_{\iota} & \\
         & F_P
     }\]
  together with Proposition \ref{prop:iota} yields
  \begin{equation*}
    \sign (\omega_i \circ \iota_0)(\nc)[(\id_A \ox (\omega_i \circ \iota_0 \circ
      \iota))^*(h_\ell)] = \sign \M[(\id_A \ox \iota)^* (h_\ell)] 
        = \sign_\iota^\M h_\ell,
  \end{equation*}
which concludes the proof.\qed
\end{proof}

\end{document}